\documentclass{conm-p-l}
\usepackage{amsmath,amssymb}

\def\char{{\rm char}}
\def\Proof{\noindent{\sl Proof.}\ }
\def\qed{{\hfill $\Box$ \medbreak}}

\newtheorem{defi}{Definition}[section]

\newtheorem{thm}[defi]{Theorem}
\newtheorem{cor}[defi]{Corollary}

\newtheorem{prop}[defi]{Proposition}
\newtheorem{example}[defi]{Example}

\DeclareMathOperator{\rad}{rad}
\DeclareMathOperator{\dl}{dl}
\DeclareMathOperator{\la}{\langle}
\DeclareMathOperator{\ra}{\rangle}

\DeclareMathOperator{\cl}{cl}
\DeclareMathOperator{\Lie}{{\mathfrak L}}

\DeclareMathOperator{\F}{\mathbb{F}}

\begin{document}

\title[Lie Properties of Restricted Enveloping Algebras]{Lie Properties of Restricted Enveloping Algebras}
\author{\textsc{Salvatore Siciliano}}
\address{Dipartimento di Matematica e Fisica ``Ennio De Giorgi", Universit\`{a} del Salento,
Via Provinciale Lecce--Arnesano, 73100--Lecce, Italy}
\email{salvatore.siciliano@unisalento.it}

\author{\textsc{Hamid Usefi}}
\address{Department of Mathematics and Statistics,
Memorial University of Newfoundland,
St. John's, NL,
Canada, 
A1C 5S7}
\email{usefi@mun.ca}

\thanks{The research of the second author was supported by NSERC of Canada. The manuscript was written while the first author visited the Department of Mathematics and Statistics at the Memorial University of Newfoundland and he expresses his gratitude for the warm hospitality during that period.}

\begin{abstract}
Let $L$ be a restricted Lie algebra over a field of positive characteristic. 
We  survey the known results about the Lie structure of the restricted enveloping algebra  $u(L)$ of $L$. Related results about the structure of the group of units
and the symmetric and skew-symmetric  elements of $u(L)$ are also discussed.
Moreover, a new theorem  about an upper bound for the Lie nilpotency class of $u(L)$ is proved.
\end{abstract}
\date{\today}

\dedicatory{ To Professor Helmut Strade on the occasion of his $70^{\text{th}}$ 
birthday.}


\maketitle

\section{Introduction}
The study of enveloping algebras that satisfy a polynomial identity (PI) was started by Latys\v{e}v in \cite{L63} in 1963 by proving that over a field of characteristic zero, the universal enveloping algebra of a Lie algebra $L$    satisfies a PI if and only if $L$ is abelian. Subsequently, Bahturin in \cite{B} completed  the characterization by dealing with the positive characteristic. Some years later, the characterization of PI restricted enveloping algebras was carried out by Passman in \cite{Pa}   and, independently, by Petrogradsky in \cite{Pe}. Since then there has been an extensive activity in this area and it is the goal of this paper to survey the main results on restricted Lie algebras satisfying certain Lie  identities.

Let $L$ be a restricted Lie algebra over a field of positive characteristic $p$ and denote by $u(L)$ the restricted (universal) enveloping algebra of $L$.
In Section \ref{Lie structure}, we exploit conditions under which $u(L)$ is Lie solvable, Lie nilpotent, and bounded Lie Engel.  The structure of the group $u(L)^\times$ of units of $u(L)$ under the assumption that $u(L)$ is an algebraic algebra are studied in  Section \ref{units}. Furthermore,  the relations between   group theoretical properties of $u(L)^\times$ and  Lie theoretical properties of $u(L)$ are exposed. It is shown that under some conditions, $u(L)$ is Lie solvable (respectively, Lie nilpotent, bounded Lie Engel) if and only if $u(L)^\times$ is solvable (nilpotent,  bounded Engel).

In Section \ref{Lie nilpotency class}, we collect the known results about the Lie nilpotency class of $u(L)$ when it is Lie  nilpotent. There, we prove a  general version of Theorem 1 of  \cite{SSp} about the bounds on the Lie nilpotency class of $u(L)$ by removing the assumptions that $L$ is finite-dimensioal and $p$-nilpotent. In Section \ref{Lie derived length}, we summarize results concerning the Lie derived length of a Lie solvable restricted enveloping algebra.

In Section \ref{involution}, we consider associative algebras with an involution. Let $A$ be an associative algebra with involution $\ast$ over a field $\F$. We denote by $A^+=\{x\in A\vert\,x^\ast=x\}$ the set of symmetric elements of $A$ and by $A^-=\{x\in A\vert\,x^\ast=-x\}$ the set of skew-symmetric elements under $\ast$.  A  question of general interest is which properties of $A^+$ or $A^-$ can be lifted to the whole algebra $A$. We consider restricted enveloping algebras endowed with the principal involution. In odd characteristic, the conditions under which  $u(L)^-$ or  $u(L)^+$ is Lie solvable, Lie nilpotent or  bounded Lie Engel are provided.  Finally, in the last section we discuss the Lie structure of ordinary enveloping algebras.

\section{Notation and definitions}\label{notation}

Let $A$ be a unital associative algebra over a field $\F$. Then   
$A$ can be regarded as a Lie algebra via the Lie bracket defined by $[x,y]=xy-yx$, for all $x,y \in A$. 
Longer Lie products in $A$ are interpreted using the left-normed convention. 
For subspaces $C,D\subseteq A$, we denote by $[C, D]$ the linear span of all elements $[c,d]$, with $c\in C$ and $d\in D$.   We say that $A$ is {\it Lie nilpotent} when $A$ is nilpotent as a Lie algebra.
If $A$ is Lie nilpotent, the {\it Lie nilpotency class} of $A$ will be denoted by $\cl_{Lie}(A)$. 
The algebra $A$ is {\it bounded Lie Engel} if there exists an integer $n$ such that $A$ satisfies the identity $[x,_n y]=[x, y,  \ldots, y]=0$, where $y$ appears $n$ times in the expression.

The $n$th {\it upper Lie power} of  $A$ is 
the ideal defined inductively by $A^{(1)}=A$ and $A^{(i)}=[A^{(i-1)},A]A$.
We say that $A$ is {\it strongly Lie nilpotent} if $A^{(i)}=0$, for some $i$, in which case 
the least index $n$ such that $A^{(n+1)}=0$  is called the strong Lie nilpotency class of $A$ and is denoted by $\cl^{Lie}(A)$.

The Lie derived series of $A$ is defined inductively by $\delta^{[0]}(A)=A$ and  $\delta^{[n+1]}(A)=[\delta^{[n]}(A),\delta^{[n]}(A)]$.  
Moreover, let us consider the series of associative ideals of $A$ defined by  $\delta^{(0)}(A)=A$  and $\delta^{(n+1)}(A)=[\delta^{(n)}(A),\delta^{(n)}(A)]A$. The algebra $A$ is said to be \emph{Lie solvable} (respectively, \emph{strongly Lie solvable}) if $\delta^{[n]}(A)=0$ ($\delta^{(n)}(A)=0$) for some $n$. In this case,  the minimal $n$ with such a property is called the $\emph{Lie derived length}$ (\emph{strong Lie derived length}) of $A$ and denoted by $\dl_{Lie}(A)$ ($\dl^{Lie}(A)$).
 The algebra $A$ is called \emph{Lie metabelian} if $\delta^{[2]}(A)=0$, and \emph{Lie center-by-metabelian} if $[\delta^{[2]}(A),A]=0$.

Clearly, strong Lie solvability implies  Lie solvability of $A$ (and $\dl_{Lie}(A)\leq \dl^{Lie}(A)$), but the converse is not true in general. For example,    one can see that the algebra $M_2(\F)$ of 2 by 2 matrices over $\F$
is Lie solvable but not strongly Lie solvable when $\F$ has characteristic 2.

Let $L$ be a restricted Lie algebra over a field $\F$ of characteristic $p>0$ and let $u(L)$ denote the restricted (universal) enveloping algebra of $L$. The terms of the lower  central series of $L$ are  $\gamma_1(L)=L$ and 
$\gamma_{n+1}(L)=[\gamma_n(L), L]$, for every $n\geq 2$. We write $L'$ for $\gamma_2(L)$.
The upper central series of $L$ is defined by $\zeta_1(L)=Z(L)$, the centre of $L$, and  
$\zeta_{n+1}(L)/\zeta_n(L)=Z(L/\zeta_n(L))$, for every $n\geq 2$.

For a subset $S$ of $L$, we denote by $\la S\ra_p$
the restricted subalgebra of $L$ generated by $S$ and by $\la S\ra_{\F}$ the subspace spanned by $S$.  We write $L'_p$ for $\la L'\ra_p$. Also, we denote by $S^{[p]^n}$ the restricted  subalgebra generated by all $x^{[p]^n}$, where $x\in S$. 
An element $x\in L$ is called \emph{$p$-nilpotent} if there exists some non-negative integer $t$ such that $x^{[p]^t}=0$; the \textit{exponent} of $x$, denoted by $e(x)$, is the least integer $s$ such that $x^{[p]^{s}}=0$.
Recall that $S$ is called \emph{$p$-nilpotent}  if there exists an integer $n$ such that $S^{[p]^n}=0$.
An element $x\in L$ is called \textit{$p$-algebraic} if $\la x\ra_p$ is finite-dimensional and  
\textit{$p$-transcendental} otherwise.

\section{The Lie structure of  restricted enveloping algebras}\label{Lie structure}

The characterization of restricted enveloping algebras has been obtained by Passman in \cite{Pa}   and, independently, by Petrogradsky in \cite{Pe}. Their result can be stated as follows:

\begin{thm}\label{passman}
 Let $L$ be a restricted Lie algebra over a field of characteristic $p>0$.  Then the restricted enveloping algebra $u(L)$ satisfies a polynomial identity if and only if $L$ has restricted subalgebras $B\subseteq A$ such that:
\begin{enumerate}
\item [{\normalfont (i)}] $\dim L/A<\infty$ and $\dim B< \infty$;
\item [{\normalfont (ii)}] $A/B$ is abelian and $B$ is central in $A$;
\item [{\normalfont (iii)}] $B$ is  $p$-nilpotent.
\end{enumerate} 
\end{thm}

The conditions under which $u(L)$ is Lie nilpotent, bounded Lie Engel, or Lie solvable are given in the following theorems.
 
\begin{thm}[\cite{RS1}]\label{Lienilpotent}
Let $L$ be a restricted Lie algebra over a field of characteristic $p>0$. The following statements are equivalent:
\begin{enumerate}
\item  $u(L)$ is Lie nilpotent;
\item $u(L)$ is strongly Lie nilpotent;
 \item  $L$ is nilpotent and $L^\prime$ is finite-dimensional and $p$-nilpotent.
\end{enumerate}
\end{thm}

\begin{thm}[\cite{RS1}]
Let $L$ be a restricted Lie algebra over a field of characteristic $p>0$. 
 Then $u(L)$ is bounded Lie Engel if and only if $L$ is nilpotent, $L^\prime$ is $p$-nilpotent, and $L$ contains a restricted ideal $I$ such that $L/I$ and $I^\prime$ are finite-dimensional. 
\end{thm}

\begin{thm}\label{thm3}
Let $L$ be a restricted Lie algebra over a field of characteristic $p>2$. The following statements are equivalent:
\begin{enumerate}
\item $u(L)$ is Lie solvable;
\item $u(L)$ is strongly Lie solvable;
\item $L^\prime$ is finite-dimensional and $p$-nilpotent.
\end{enumerate}
\end{thm}

The equivalence of (1) and (3) in Theorem \ref{thm3} is shown  in \cite{RS1} whereas it is shown in \cite{S1} that (2) and (3) are equivalent for all $p>0$.  
On the other hand,  the characterization of Lie solvable restricted enveloping algebras in characteristic 2 was carried out only  recently. The main difficulty here is that unlike other characteristics, Lie solvability is not a non-matrix polynomial identity in characteristic 2.
Recall that a polynomial identity is called non-matrix  if it is not satisfied by the algebra $M_2(\F)$. 
Indeed, if $\char \F=2$  then $M_2(\F)$ is Lie center-by-metabelian.
Using the standard PI-theory, like Posner's Theorem, one can deduce that if $R$ is an associative algebra that satisfies a non-matrix PI over a field $\F$ of characteristic $p$   then $[R, R]R$ is nil. If we further assume that $R$ is Lie solvable and $p\neq 2$,  then $[R, R]R$ is nil of bounded index (see \cite{R97}).
Moreover, if we restrict ourselves to $R=u(L)$ then  $R$ satisfies a non-matrix PI if and only if
$[R, R]R$ is nil of bounded index (see \cite{RW99}). However,   if $u(L)$ is Lie solvable and $p=2$ then 
$L^\prime$ may not be even  nil as the following theorem states. We recall that a restricted Lie algebra is said to be {\it strongly abelian} if it is abelian and its power mapping is trivial.

\begin{thm}[\cite{SU}]\label{solvchar2}
Let  $L$ be a   restricted Lie algebra over a field $\F$ of characteristic 2. Let $\bar{\F}$ be the algebraic closure of $\F$ and set $\Lie=L \otimes_{\F} \bar{\F}$.  Then $u(L)$ is Lie solvable if and only if  
$\Lie$ has  a finite-dimensional 2-nilpotent restricted ideal $I$ such that $\bar \Lie=\Lie/I$ satisfies one of the following conditions:
\begin{enumerate} 
\item[{\normalfont (i)}]  $\bar \Lie$ has an abelian restricted ideal of codimension at most $1$;
\item[{\normalfont (ii)}]  $\bar \Lie$ is nilpotent of class 2 and $\dim \bar \Lie/Z(\bar \Lie)=3$;
\item[{\normalfont (iii)}]  $\bar \Lie=  \la  x_1, x_2, y\ra_{\bar{\F}} \oplus Z(\bar \Lie)$, where  $[x_1,y]=x_1$,  $[x_2,y]=x_2$, and  $[x_1,x_2]\in Z(\bar \Lie)$;
\item [{\normalfont (iv)}]  $\bar \Lie=  \la  x, y\ra_{\bar{\F}} \oplus H \oplus Z(\bar \Lie)$, where  $H$ is a strongly abelian finite-dimensional restricted subalgebra of $\bar \Lie$ such that $[x,y]=x$, $[y,h]=h$, and $[x,h]\in Z(\bar \Lie)$  for every $h\in H$;
\item [{\normalfont (v)}] $\bar \Lie= \la  x, y\ra_{\bar{\F}} \oplus H \oplus Z(\bar \Lie) $, where $H$ is a finite-dimensional abelian subalgebra of $\bar \Lie$ such that $[x,y]=x$, $[y,h]=h$, $[x,h]\in Z(\bar \Lie)$, and $[x,h]^{[2]}=h^{[2]}$, for every $h\in H$.
\end{enumerate}
\end{thm}

Note that the cases (ii)-(v) can occur only when $L^\prime$ is finite-dimensional. In other words, if $u(L)$ is Lie solvable and $L^\prime$ is infinite-dimensional, then $L$ has a 2-abelian restricted ideal of codimension at most 1. The  following example shows that   the extension of the ground field is necessary in Theorem \ref{solvchar2}.

\begin{example}[\cite{SU}]\label{examplecod} \emph{ Let $\F$ be a field of characteristic 2 containing two elements  $\alpha, \beta$ such that the following condition holds: If $\lambda_1,\lambda_2,\lambda_3$ are  in $\F$ and  $\lambda_1^2+\lambda_2^2 \alpha + \lambda_3^2 \beta=0$ then $\lambda_1=\lambda_2=\lambda_3=0$. For instance, one can consider the field $\mathbb{K}(X,Y)$ of rational  functions in two indeterminates over any field $\mathbb{K}$ of characteristic 2, and $\alpha=X$ and $\beta=Y$.  
Let $L$ be the $\F$-vector space  having the elements $x, x_1, x_2, x_3, z_1, z_2, z_3$ as basis. We define a 
 restricted Lie algebra structure on $L$ by setting  
$[x,x_1]=[x, x_3]=z_1$, $[x,x_2]=z_2$,  $ [x_1,x_2]=z_3$, $[x_1, x_3]=\frac{\beta}{\alpha}z_3$, $[x_2,x_3]=0$, $ z_1^{[2]}= z_1, z_2^{[2]}= \alpha  z_1$, $ z_3^{[2]}=\beta z_1$ and $z_i \in Z(L)$, $x^{[2]}=x_i^{[2]}= 0$, for $i = 1, 2, 3$. Then $\Lie=L \otimes_{\F} \bar{\F}$ contains a 2-nilpotent restricted ideal $J$ such that $\Lie/J$ has an abelian restricted ideal of codimension 1, thus $u(\Lie)$ is Lie solvable. On the other hand, $L$ does not 
contain any  restricted ideal $I$ such that $L/I$ satisfies one of the five conditions of Theorem \ref{solvchar2}.}
\end{example}

It should be mentioned that these problems have been investigated in the more general setting of 
restricted Lie superalgebras, as well. For instance, the characterization of restricted Lie superalgebras whose enveloping algebras satisfy a PI was carried out in \cite{P92}. Furthermore, restricted Lie superalgebras whose enveloping algebras satisfy a non-matrix PI have been recently characterized in \cite{UJPAA, UJA}.

\section{The group of units}\label{units}

Let $A$ be a unital associative algebra over a
field $\F$. We shall denote by $A^\times$ the
group of units of $A$. Unlike group algebras very few results are known about the group of units of restricted 
enveloping algebras. For example, it is not even known when $u(L)^\times$ is abelian in general. However, if we restrict ourselves to the case that $u(L)$ is algebraic then we can state the following results.

\begin{thm}[\cite{JRS}]
Let $L$ be a restricted Lie algebra over a field of
odd characteristic $p$ with at least 5 elements. If
$u(L)$ is algebraic then the following conditions
are equivalent.
\begin{enumerate}
\item The group of units $u(L)^\times$ is solvable.
\item The algebra $u(L)$ is Lie solvable.
\item
The derived subalgebra $[L,L]$ is both
finite-dimensional and $p$-nilpotent.
\end{enumerate}
\end{thm}

\begin{thm}[\cite{JRS}]\label{units-engel}
Let $L$ be a restricted Lie algebra over a perfect field
 of positive characteristic $p$ with at least 3
elements. If $u(L)$ is algebraic then the following
conditions are equivalent.
\begin{enumerate}
\item The group of units $u(L)^\times$ is bounded Engel.
\item The algebra $u(L)$ is bounded Engel.
\item
The restricted Lie algebra $L$ is nilpotent, $L$
contains a restricted ideal $I$ such that $L/I$ and
$[I,I]$ are finite-dimensional, and $[L,L]$ is
$p$-nil of bounded index.
\end{enumerate}
\end{thm}

\begin{thm}[\cite{JRS}]\label{units-nilpotent}
Let $L$ be a restricted Lie algebra over a perfect field of
positive characteristic $p$ with at least 3
elements. If $u(L)$ is algebraic then the following
conditions are equivalent.
\begin{enumerate}
\item The group of units $u(L)^\times$ is nilpotent.
\item The algebra $u(L)$ is Lie nilpotent.
\item
The restricted Lie algebra $L$ is nilpotent and
$[L,L]$ is both finite-dimensional and
$p$-nilpotent.
\end{enumerate}
\end{thm}

The ground field $\F_2$ was correctly omitted in
Theorems \ref{units-engel} and \ref{units-nilpotent}. Indeed, as it is pointed out in \cite{JRS}, one can consider the restricted
enveloping algebra $u(L)$ of the restricted Lie
algebra $L$ over $\F_2$ with a basis $\{x,y\}$ such
that $[x,y]=x$, $x^{[2]}=0$, and $y^{[2]}=y$. Then
$u(L)^\times$ is isomorphic to the Klein four group.
Thus, $u(L)^\times$ is abelian even though $u(L)$ is
not bounded Lie Engel.

\section{Lie nilpotency class}\label{Lie nilpotency class}
Following \cite{RS2}, for every positive integer $m$ 
we consider the restricted ideal of $L$ given by
$$
D_{(m)}(L)=L\cap u(L)^{(m)}.
$$
The terms of this sequence are called \emph{the upper Lie 
dimension subalgebras} of $L$.
In analogy to modular group rings  (cf. \cite{PS}),  
the $D_{(m)}(L)$'s can be described as  follows
(see \cite[Theorem 4.4]{RS2}):
\begin{equation}\label{shalev}
D_{(m+1)}(u(L))=\left \{ \begin{array}{ll}
L &  \quad  m=0\\\nonumber
L^{\prime}_{p}& \quad  m=1\\ 
{D_{(\lceil \frac{m+p}{p}\rceil)}(L)}^{[p]}+[D_{(m)}(L),L] & \quad m\geq 2
\end{array} \right.
\end{equation}
where, for a real number $r$, $\lceil r \rceil$ denotes the smallest integer not less than $r$.
Moreover, it holds that 
\begin{align}\label{Dm}
D_{(m+1)}(L)=\sum_{(i-1)p^{j}\geq m}{{\gamma}_{i}(L)}^{[p]^{j}}.
\end{align}

For every $m\geq 1$, put $d_{(m)}=\dim_{\F}(D_{(m)}(L)/D_{(m+1)}(L))$. Then, we have the following:

\begin{thm}[\cite{RS2}]\label{upper-nilpotent-class} Let $L$ be a restricted Lie algebra over a field of characteristic $p>0$. If $u(L)$ is strongly  Lie nilpotent then:
$$
\cl^{Lie}(u(L))=1+(p-1) \sum_{m\geq 1} md_{(m+1)}.
$$
\end{thm}

Moreover, as the following theorem shows, the two Lie nilpotency classes coincide when $p>3$.
\begin{thm}[\cite{RS2}]\label{RS2}  Let $L$ be a restricted Lie algebra of characteristic
$p>3$ with $u(L)$ being Lie nilpotent.  
Then  $\cl_{Lie}(u(L))=\cl^{Lie}(u(L))$.
\end{thm}

It is still unknown whether Theorem \ref{RS2} holds true in  characteristics $p=2, 3$.
Under the assumption that $L$ is finite-dimensional and $p$-nilpotent,  the authors in \cite{SSp} gave a lower bound for $\cl_{Lie}(u(L))$ and an upper bound for $\cl^{Lie}(u(L))$. Now we can prove their result in general.

\begin{thm}\label{nilpclass}
Let $L$ be a
restricted Lie algebra over  a field of characteristic $p>0$. If $u(L)$ is Lie nilpotent then  
$$
p^{e([x,y])}\leq \cl_{Lie}(u(L)) \leq
\cl^{Lie}(u(L)) 
\leq p^{\dim_{\F} L^\prime_p},
$$
for every $x,y\in L$.
\end{thm}
\Proof By Proposition 1 in \cite{SSp}, we know that 
$p^{e([x,y])}\leq \cl_{Lie}(u(L))$ for every restricted Lie algebra $L$ and every  $x,y\in L$.
Now we prove that $\cl^{Lie}(u(L)) 
\leq p^{\dim_{\F} L^\prime_p}$ for every restricted Lie algebra $L$.

Without loss of generality we can assume that the ground field $\F$ is algebraically closed. We first assume that $L$ is finite-dimensional. As $L$ is nilpotent the semisimple elements of $L$ are central and they form the unique maximal torus $T$, which  is also a restricted ideal of $L$.  Denote by $\rad_p(L)$ the $p$-radical of $L$.   Let $I$ denote the set of all $p$-nilpotent elements of $L$. Let $x,y\in I$. For every positive integer $k$  one has that 
$$(x+y)^{[p]^k}\equiv x^{[p]^k}+y^{[p]^k} \mod L^\prime_p,$$ 
so that  $(x+y)^{[p]^n}\in L^\prime_p$ for a sufficiently large $n$. Since $L^\prime_p$ is $p$-nilpotent
 by Theorem \ref{Lienilpotent}, this forces that $x+y$ is $p$-nilpotent and $I$ is a subspace of $L$.  As $L^\prime_p\subseteq I$, we have that $I$ is a $p$-nilpotent restricted ideal of $L$, which entails that $\rad_p(L)=I$. Now, as the ground field is algebraically closed, for every element $x$ of $L$ we can consider the Jordan-Chevalley decomposition $x=x_s+x_n$ (see \cite[\S 2, Theorem 3.5]{SF}), where $x_s$ and $x_n$ are the semisimple  and $p$-nilpotent parts of $x$, respectively. This proves that $L=\rad_p(L) \oplus T$ and so $u(L)\cong u(T)\otimes_{\F}u(\rad_p(L))$.  Note that, as $T$ is central, for every $x,y\in L$ we have $[x,y]=[x_n,y_n]$ and so $[L,L]=[\rad_p(L),\rad_p(L)]$. As a consequence, by Theorem 1 of \cite{SSp}, for every $x,y\in L$ we have
\begin{align*}
\cl^{Lie}(u(L))
= \cl^{Lie}(u(\rad_p(L)))\leq p^{\dim_{\F} L^\prime_p},
\end{align*}
and the claim follows for the finite-dimensional case. 

Now suppose that $L$ is an arbitrary restricted Lie algebra. We claim that there exists a finite-dimensional nilpotent restricted Lie algebra $\mathcal{L} $ such that $\mathcal{L} ^\prime_p$  and $L^\prime_p$ are isomorphic as restricted Lie algebras. Suppose that the claim holds. 
First observe by Equation \eqref{Dm} that $D_{(m)}(L)\cong D_{(m)}(\mathcal{L} )$.
Then, by Theorem \ref{upper-nilpotent-class}, we have
\begin{align*}
\cl^{Lie}(u(L)) &=1+(p-1) \sum_{m\geq 1} m \dim_{\F}(D_{(m)}(L)/D_{(m+1)}(L))\\
&=1+(p-1) \sum_{m\geq 1} m \dim_{\F}(D_{(m)}(\mathcal{L} )/D_{(m+1)}(\mathcal{L} ))\\
&=\cl^{Lie}(u(\mathcal{L} )).
\end{align*}
Since $\mathcal{L} $ is finite-dimensional,  we deduce by the first part of the proof that 
$$
 \cl^{Lie}(u(L)) =\cl^{Lie}(u(\mathcal{L} )) \leq  p^{\dim_{\F} \mathcal{L}^\prime_p}=p^{\dim_{\F} L^\prime_p}.
$$

Thus, it is enough to prove the claim.
As $L^\prime$ is finite-dimensional there exist $a_1,\ldots,a_n,b_1,\ldots,b_n \in L$ such that $L^\prime$ is spanned by the commutators $[a_i,b_i]$. Now consider the restricted subalgebra $H$ of $L$ generated by the $a_i,b_j$, $i,j=1,\ldots,n$. Note that $H^\prime_p=L^\prime_p$ and by well-known results about the structure of finitely generated abelian restricted Lie algebras over perfect fields (see e.g. Section 4.3 in \cite{BMPZ}), we deduce that there exists $p$-algebraic elements $x_1, \ldots, x_r$ and $p$-transcendental elements $y_1,\ldots, y_s$ in $H$ such that  
$$
H/H^\prime_p \cong \la \bar x_1 \ra_p \oplus \cdots \oplus \la \bar x_r \ra_p \oplus \la \bar y_1 \ra_p \oplus \cdots \oplus \la \bar y_s \ra_p,
$$ 
where $\bar x_i$ and $\bar y_j$ are, respectively, the images of $x_i$ and $y_j$ in  $H/H^\prime_p$.
Let $\{z_1,\ldots,z_t\}$ be an $\F$-basis for $H^\prime_p$. Now let $e$ be the minimal integer such that $p^e$ is greater than the nilpotency class of $L$. Theorem 2.3 and Proposition 2.1 in \cite{SF} allow us  to  define a restricted Lie algebra $\mathcal{L}$ whose underlying Lie algebra is $H$ and its $p$-map which we denote  by $[p]_1$ satisfies the following:
\begin{align*}
x_i^{[p]_1^m}=x_i^{[p]^m},\quad z_j^{[p]_1}=z_j^{[p]},
\end{align*}
for all $i=1, \ldots, r, j=1,\ldots, t$,  all positive integers $m$,
and 
$$
y_k^{[p]_1^m}=
\begin{cases}
y_k^{[p]^m} , &\text{if $m<e$};\\
0 , &\text{if $m\geq e$},\\
\end{cases}
$$
for all $k=1, \ldots, s$. We observe that $\mathcal{L}$ is finite-dimensional and $\mathcal{L}^\prime_p\cong H'_p$, as restricted Lie algebras, which completes the proof.
\qed

As a consequence, one has that $\cl_{Lie}(u(L))=\cl^{Lie}(u(L))=p^{\dim_{\F} L^\prime_p}$  
when $L^\prime_p$ is cyclic. Conversely, it follows from \cite{SSp} that if the nilpotency class of $L$ is not greater than $p$, then $\cl_{Lie}(u(L))$ reaches the maximal value $p^{{\dim}_{\F} L^\prime_p}$ if and only if $L^\prime_p$ is cyclic. This is no longer true in general without the assumption on the nilpotency class of $L$ (see Example 1 in \cite{SSp}).

\section{Lie derived length}\label{Lie derived length}
Determining the Lie derived lengths of a Lie solvable restricted enveloping algebra is in general a difficult task. In this sections we summarize the known results in this direction.
An associative algebra $A$ is called \emph{Lie metabelian} (respectively, \emph{strongly Lie metabelian}) if $\delta^{[2]}(A)=0$ ($\delta^{(2)}(A)=0$), and \emph{Lie center-by-metabelian} if $[\delta^{[2]}(A),A]=0$.
Lie metabelian restricted enveloping algebras were characterized in \cite{S2}: 

\begin{thm}\label{Lie metabelian}
Let $L$ be a non-abelian restricted Lie algebra over a field $\F$
of characteristic $p>0$. Then the following three conditions are equivalent:
\begin{enumerate}
\item [{\normalfont 1)}] $u(L)$ is Lie metabelian; 
\item [{\normalfont 2)}] $u(L)$ is strongly  Lie metabelian; 
\item [{\normalfont 3)}]  one of the following conditions is satisfied:
 {\begin{enumerate}
\item [{\normalfont (i)}] $p=3$, $\dim_{\F} L^\prime =1$,  $L^\prime$
is central  and ${L^\prime}^{[p]}=0$; \item [{\normalfont (ii)}]
$p=2$, $\dim_{\F} L^\prime = 2$, $L^\prime$ is central and
${L^\prime}^{[p]}=0$; \item [{\normalfont (iii)}] $p=2$, $\dim_{\F} L^\prime =1$   and ${L^\prime}^{[p]}=0$
\end{enumerate}}
\end{enumerate}
\end{thm}

 Moreover, for $p>3$,  it is shown in \cite{RT}  that $u(L)$ is Lie center-by-metabelian if and only if $L$ is abelian. In \cite{S2}  the characterization of Lie center-by-metabelian in odd characteristic has been completed by settling the more difficult case $p=3$. We have the following:
 
 \begin{thm}[\cite{S2}]\label{Liecentrm}
Let $L$ be a restricted Lie algebra over a field $\F$ of characteristic $p>2$. 
Then $u(L)$ is Lie center-by-metabelian if and only if either $L$ is abelian or $p=3$, $\dim_{\F} L^\prime=1$, $L^\prime$ is central, and ${L^\prime}^{[p]}=0$. 
\end{thm}

As a consequence of Theorem \ref{Liecentrm},  in odd characteristic a Lie center-by-metabelian restricted enveloping algebra is in fact Lie metabelian.  This is no longer true for $p=2$, for instance consider $L=\la x, y, z\mid [x,y]=z, [x,z]=[y,z]=x^{[2]}=y^{[2]}=0, z^{[2]}=z\rangle$.

If $u(L)$ is strongly Lie solvable then an upper bound for the strong Lie derived length of $u(L)$ is provided by the following result:

\begin{prop}[\cite{S1}]\label{log}
Let $L$ be a restricted Lie algebra over a field of characteristic
$p>0$. If $u(L)$ is strongly Lie solvable then
$\dl^{Lie}(u(L))\leq \lceil \log_2 (2t(L_p^\prime))\rceil$.
\end{prop}
Here, $t(L_p^\prime)$ is the nilpotency index of the augmentation ideal of $u(L_p^\prime)$, which can be computed by  Corollary 2.4 of   \cite{RS2}. Furthermore,  the following theorem provides a lower bound for the Lie derived length of non-commutative restricted enveloping algebras.
  
\begin{thm}[\cite{S1}]
Let $L$ be a non-abelian restricted Lie algebra over a field of
characteristic $p>0$.  If $u(L)$ is Lie solvable then $\dl_{Lie}(u(L))\geq \lceil
\log_2(p+1)\rceil$.
\end{thm}

A characterization of non-commutative restricted enveloping algebras whose Lie derived length coincides with the above minimal value has been subsequently carried out in \cite{CSS}:

\begin{thm}\label{derived-length}
Let $L$ be a restricted Lie algebra over a field $\F$
of characteristic $p>0$. Then the following three conditions are equivalent:
\begin{enumerate}
\item [{\normalfont 1)}] $\dl_{Lie}(u(L))=\lceil\log_2(p+1)\rceil$; 
\item [{\normalfont 2)}] $\dl^{Lie}(u(L))=\lceil\log_2(p+1)\rceil$; 
\item [{\normalfont 3)}] one of the following conditions is satisfied:
 {\begin{enumerate}
\item [{\normalfont (i)}] $p=2$, $\dim_{\F} L^\prime = 2$, $L^\prime$ is central, and  ${L^\prime}^{[p]}=0$;
\item [{\normalfont (ii)}] $p=2$, $\dim_{\F} L^\prime=1$, and ${L^\prime}^{[p]}=0$;
\item [{\normalfont (iii)}] $p>2$, $\dim_{\F} L^\prime=1$, $L^\prime$ is central, and  ${L^\prime}^{[p]}=0$. 
\end{enumerate}}
\end{enumerate}
\end{thm}

Let $L$ be a non-abelian restricted Lie algebra over a field of characteristic $p>2$. If $u(L)$
is Lie nilpotent then, by Theorem \ref{nilpclass}, the Lie nilpotency class $\cl_{Lie}(u(L))$
of $u(L)$ is at least $p$.  From Theorem \ref{derived-length}  and Corollary 1 of \cite{S2}  we deduce     the following:
  
 \begin{cor}
 Let $L$ be a restricted Lie algebra over a field of characteristic $p>2$. Then 
 $\dl_{Lie}(u(L))=\lceil\log_2(p+1)\rceil$ if and only if $\cl_{Lie}(u(L))=p$.
 \end{cor}
 
Note that, at least in characteristic $2$, if $u(L)$ is strongly Lie solvable then it is possible to have $\dl_{Lie}(u(L))< \dl^{Lie}(u(L))$, as the following example shows.
\begin{example}\emph{
Let $L$ be the 5-dimensional restricted Lie algebra over a field $\F$ of characteristic $2$ given by
$$
L=\langle x,y,z,v,w\mid \, [x,y]=z, z,v,w\in Z(L), x^{[p]}=y^{[p]}=w^{[p]}=0,  z^{[p]}=v, v^{[p]}=w \rangle. 
$$
Then one can see that
$3=\dl_{Lie}(u(L))< \dl^{Lie}(u(L))=4$.
}
\end{example}

However, we do not know any example of a  restricted enveloping algebra defined over a field of odd characteristic whose Lie derived lengths are different.

Some other results on the derived length of $u(L)$ under the assumption that $L$ is powerful  can be found in \cite{SW}.

\section{Skew-symmetric and symmetric elements}\label{involution}

Let $A$ be an algebra with involution $\ast$ over a field $\F$. We denote by $A^+=\{x\in A\vert\,x^\ast=x\}$ the set of symmetric elements of $A$ under $\ast$ and by $A^-=\{x\in A\vert\,x^\ast=-x\}$ the set of skew-symmetric elements.  A  question of general interest is which properties of $A^+$ or $A^-$ can be lifted to the whole algebra $A$.  Herstein conjectured in \cite{H67, H76} that if  the symmetric or skew-symmetric elements of  a ring $R$ with involution satisfy a polynomial identity, then so does $R$. This conjecture was proved  by  Amitsur in  \cite{A68} and subsequently generalized by him in  \cite{A69}:

\begin{thm}[\cite{A68}]\label{amitsur} 
Let $A$ be an associative algebra with involution. If $A^+$ (or $A^-$) satisfies a 
polynomial identity then $A$ satisfies a polynomial identity.
\end{thm}

 A subset $S$ of $A$ is said to be Lie nilpotent if there exists a positive integer $n$ such that 
$$[x_1,\ldots,x_n]=0$$
 for every $x_1,\ldots,x_n\in S$, while  $S$ is said to be bounded Lie Engel if there is an $n$ such that 
 $$[x,_ny]=0$$
  for every $x,y\in S$. Also,  we put $[x_1,x_2]^o =[x_1,x_2]$ and 
$$
[x_1,x_2,\ldots,x_{2^{n+1}}]^o=[[x_1,\ldots,x_{2^n}]^o,
[x_{2^n+1},\ldots,x_{2^{n+1}}]^o].
$$ 
The subset $S$  is said to be Lie solvable if there exists an $n$ such that 
$$
[x_1,x_2,\ldots,x_{2^{n+1}}]^o=0,
$$
 for every $x_1,\ldots,x_{2^{n+1}}\in S$.

Now consider the group algebra $\F G$ of a group $G$ under the canonical involution induced by $g \mapsto g^{-1}$, for every $g\in G$. There is an extensive literature devoted  to establish the extent to which the symmetric or skew-symmetric elements of  $\F G$ under  a given involution  determine the Lie structure of the group algebra, for example see \cite{CLS, GS93, GS06, L99, L00, LSS, LSS1}.

Let $L$ be a restricted Lie algebra  over a field ${\mathbb F}$ of characteristic $p>2$ and let $u(L)$ be the restricted enveloping algebra of $L$.  
We denote by $\top$  the \emph{principal involution} of $u(L)$, that is, the unique ${\mathbb F}$-antiautomorphism of $u(L)$ such that $x^\top=-x$ for every $x$ in $L$. We recall that $\top$ is just the antipode of the  ${\mathbb F}$-Hopf algebra $u(L)$.

The conditions under which  $u(L)^-$ or  $u(L)^+$ is Lie solvable, Lie nilpotent or  bounded Lie Engel are determined in \cite{Skew, SUIsr}. These  results are summarized in the following three theorems.

\begin{thm}\label{bengel}
Let $L$ be a restricted Lie algebra over a field ${\mathbb F}$ of characteristic $p>2$. Then the following conditions are equivalent:
\begin{enumerate}
\item [{\normalfont 1)}] $u(L)^+$ is bounded Lie Engel;
\item [{\normalfont 2)}] $u(L)^-$ is bounded Lie Engel;
\item [{\normalfont 3)}] $u(L)$  is bounded Lie Engel;
\item [{\normalfont 4)}] $L$ is nilpotent, $L^\prime$ is $p$-nilpotent, and $L$
contains a restricted ideal $I$ such that $L/I$ and $I^\prime$ are finite-dimensional.
\end{enumerate}
\end{thm}

\begin{thm}
Let $L$ be a restricted Lie algebra over a field ${\mathbb F}$ of characteristic $p>2$. Then the following conditions are equivalent:
\begin{enumerate}
\item [{\normalfont 1)}] $u(L)^+$ is Lie nilpotent;
\item [{\normalfont 2)}] $u(L)^-$ is Lie nilpotent;
\item [{\normalfont 3)}] $u(L)$ is Lie nilpotent;
\item [{\normalfont 4)}] $L$ is nilpotent and $L^\prime$ is finite-dimensional and $p$-nilpotent.
\end{enumerate}
\end{thm}

\begin{thm}\label{solvable}
Let $L$ be a restricted Lie algebra over a field ${\mathbb F}$ of characteristic $p>2$. Then the following conditions are equivalent:
\begin{enumerate}
\item [{\normalfont 1)}] $u(L)^+$ is Lie solvable;
\item [{\normalfont 2)}] $u(L)^-$ is Lie solvable;
\item [{\normalfont 3)}] $u(L)$ is Lie solvable;
\item [{\normalfont 4)}] $L^\prime$ is finite-dimensional and $p$-nilpotent.
\end{enumerate}

\end{thm}

\section{Ordinary enveloping algebras}
Let $L$ be a Lie algebra over an arbitrary field $\F$ and denote by $U(L)$ the universal enveloping algebra of $L$.
Latys\v{e}v in \cite{L63} proved that over a field of characteristic zero,  $U(L)$  is PI if and only if $L$ is abelian. Bahturin in \cite{B} extended Latyshev's result to  the positive characteristic 
by proving that $U(L)$ is PI if and only if $L$ has an abelian ideal  of finite codimension and the adjoint representation of $L$ is algebraic of bounded degree. 

Suppose that  $p>0$ and consider the universal $p$-envelope 
of $L$
$$
\hat{L}=\sum_{k \geq 0} {L^{p^k}}\subseteq U(L),
$$
where $L^{p^k}$ is the ${\mathbb F}$-vector space 
spanned by the set $\{l^{p^k}\vert\, l\in L \}$.
Then $\hat{L}$ is a restricted Lie algebra with $p$-map given by $h^{[p]}=h^p$ for all $h\in \hat{L}$. Note that by Corollary $1.1.4$ of \cite{S},  $U(L)=u(\hat{L})$. Hence, we can apply  Theorem \ref{solvchar2}, to deduce the following:

\begin{thm}[\cite{SU}] Let $L$ be a Lie algebra over a field ${\mathbb F}$ of characteristic 2. If $U(L)$ is Lie solvable then one of the following conditions is satisfied:
\begin{enumerate}
\item [{\normalfont (i)}]  ${L}$ contains an abelian ideal of codimension at most 1;
\item [{\normalfont (ii)}] $L$ is nilpotent of class 2 and $\dim_{\F}L/Z(L)=3$;
\item [{\normalfont (iii)}]  $L=\la x_1, x_2, y\ra_{\F}\oplus Z(L)$, with $[x_1,y]=x_1$, $[x_2,y]=x_2$, and $[x_1,x_2]\in Z(L)$.
\end{enumerate}
\end{thm}

Now suppose that $U(L)^+$ or $U(L)^-$ (with respect to the principal involution) is Lie solvable or bounded Lie Engel over a field of characteristic $p\neq 2$.  Then, by Theorem \ref{amitsur},  $U(L)$ satisfies a polynomial identity. Hence,  if $p=0$  then $L$ is abelian    by    Latys\v{e}v's Theorem.
On the other hand, if $p>2$ then,  by Theorems \ref{solvable} and \ref{bengel}, we have   that $(\hat{L})^\prime$ is $p$-nilpotent. Since $u(\hat{L})=U(L)$ has no nontrivial zero divisors, we conclude that $L^\prime=0$. 
This proves the following:
 
\begin{cor} Let $L$ be a Lie algebra over a field ${\mathbb F}$ of characteristic $p\neq 2$. Then the following conditions are equivalent:
\begin{enumerate}
\item $U(L)$ is Lie solvable;
\item $U(L)$ is bounded Lie Engel;
\item $U(L)^+$ or $U(L)^-$ is Lie solvable;
\item $U(L)^+$or $U(L)^-$ is bounded Lie Engel;
\item $L$ is abelian.
\end{enumerate}
\end{cor}  

We conclude by mentioning that Lie superalgebras whose enveloping algebras satisfy a PI have been described  in  \cite{B85, P92}. Furthermore, the characterization of Lie superalgebras whose enveloping algebras satisfy a non-matrix PI was recently carried out in \cite{BRU}.

\enddocument